\documentclass[11pt]{article}
\usepackage{epsfig}
\usepackage{amssymb,amsmath}
\usepackage{mathtools}
\usepackage{graphicx, epstopdf}
\usepackage{caption, floatrow}
\usepackage[top=0.9in, bottom=0.9in, left=0.85in, right=0.85in]{geometry}

\begin{document}

\title{\bf On numerical modelling of contact lines in fluid flows}
\author{Dmitry E. Pelinovsky$^{1}$ and Chengzhu Xu$^{2}$ \\
    {\em $^1$Department of Mathematics, McMaster University, Hamilton ON, Canada, L8S 4K1} \\
    {\em $^2$Department of Applied Mathematics,  University of Waterloo, Waterloo ON, Canada, N2L 3G1}}

\date{\today}
\maketitle

\begin{abstract}
We study numerically a reduced model proposed by Benilov and Vynnycky (J. Fluid Mech. {\bf 718} (2013), 481), who examined
the behavior of a contact line with a $180^{\circ}$ contact angle between liquid and a moving plate,
in the context of a two-dimensional Couette flow. The model is given by a linear fourth-order
advection-diffusion equation with an unknown velocity, which is to be determined
dynamically from an additional boundary condition at the contact line.

The main claim of Benilov and Vynnycky is that for any physically relevant initial condition,
there is a finite positive time at which the velocity of the contact line tends to negative infinity,
whereas the profile of the fluid flow remains regular.
Additionally, it is claimed that the velocity behaves as the logarithmic function of time near the blow-up time.

Compared to the previous computations based on COMSOL built-on algorithms,
we use MATLAB software package and develop a direct finite-difference method to study
dynamics of the reduced model under different initial conditions. We confirm the first claim
but also show that the blow-up behavior is better approximated by a power function,
compared with the logarithmic function. This numerical result suggests a simple explanation
of the blow-up behavior of contact lines.
\end{abstract}

\section{Introduction}

Contact lines are defined by the triple-point intersection of the rigid boundary, fluid flow
and the vacuum state. Flows with the contact line at $180^{\circ}$ contact angle were
discussed in \cite{Benney,Ngan}, where corresponding solutions of the Navier-Stokes equations
were shown to have no physical meanings. In the recent paper, Benilov and Vynnycky
\cite{Benilov} analyzed the behavior of the contact line asymptotically by using
the thin film equations.

\begin{figure}[H]
\includegraphics[height=2in]{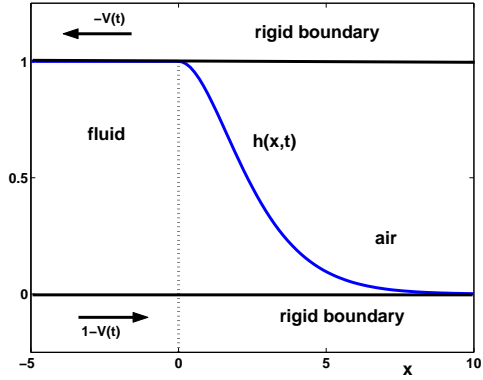}
\caption{A two-dimensional Couette flow with a free boundary, in the reference frame co-moving with the contact line.}
\label{fig-scheme}
\end{figure}

Consider a two-dimensional Couette flow shown on Figure \ref{fig-scheme}, where two horizontal rigid plates
are separated by a distance normalized to unity, with the lower plate moving to the right relatively to
the upper plate with a velocity normalized to unity. The space between the plates is filled
with an incompressible fluid on the left, and vacuum (that is, gas with negligible density) on the right,
separated by a free boundary. The $x$-axis is directed along the lower plate,
and the contact line is located on the upper plate.

Physically relevant flows correspond to the configuration, where
the fluid-filled region to the right of the contact line decays monotonically, and is carried away by the lower plate
to some residual thickness $h_{\infty}$ as $x \rightarrow \infty$. The velocity of the contact line is $V(t)$
and the reference frame on Figure \ref{fig-scheme} moves to the left with the velocity $V(t)$ so that
the contact line is placed dynamically at the point $x = 0$. Note that the velocity $V(t)$ is an
unknown variable to be found as a function of time $t$. The shape of the fluid-vacuum interface
at time $t$ is described by the graph of the function $y = h(x,t)$ for $x > 0$,
where $h$ is the thickness of fluid-filled region.

By using asymptotic analysis and the lubrication approximation, Benilov and Vynnycky
\cite{Benilov} derived the following nonlinear advection--diffusion equation
for the free boundary $h(x,t)$ of the fluid flow:
\begin{equation}
\label{model}
\frac{\partial h}{\partial t} + \frac{\partial}{\partial x} \left[
\frac{h^3}{3} \left( \alpha^3 \frac{\partial^3 h}{\partial x^3}
+ \frac{\partial h}{\partial x} \right) + (1 - V(t)) h \right] =
0, \quad x > 0, \;\; t > 0,
\end{equation}
The boundary conditions $h|_{x = 0} = 1$ and $h_x |_{x = 0} = 0$ define the normalized thickness
and the contact line location, whereas the flux conservation gives the boundary condition
for $h_{xxx} |_{x = 0} = -\frac{3}{2 \alpha^3}$. Here and henceforth, we use the subscript to
denote the partial derivative. In addition, we fix $\alpha^3 = 3$ for convenience.
Existence of weak solutions of the thin-film equation (\ref{model}) for constant values of $V$
and Neumann boundary conditions on a finite interval was recently
shown by Chugunova {\em et al.} \cite{Chugunova1,Chugunova2}.

Using further asymptotic reductions with
\begin{equation}
\label{asymptotic-reduction}
h - 1 = \mathcal{O}(|V|^{-1}), \quad x = \mathcal{O}(|V|^{-1/3}), \quad
t = \mathcal{O}(|V|^{-4/3}), \quad \mbox{\rm as} \quad |V| \to \infty,
\end{equation}
Benilov and Vynnycky \cite{Benilov} reduced the nonlinear equation (\ref{model})
with $\alpha^3 = 3$ to the linear advection--diffusion equation:
\begin{equation}
\label{pde}
\frac{\partial h}{\partial t} + \frac{\partial^4 h}{\partial x^4} =
V(t) \frac{\partial h}{\partial x}, \quad x > 0, \;\; t > 0,
\end{equation}
subject to the boundary conditions
\begin{equation}
\label{bc-pde}
h |_{x = 0} = 1, \quad  h_{x} |_{x = 0} = 0, \quad h_{xxx} |_{x = 0} = - \frac{1}{2}, \quad
t \geq 0,
\end{equation}

Physically relevant solutions corresponds to the monotonically decreasing
solutions with  $h \to h_{\infty}$ and $h_x, h_{xx} \to 0$ as $x \to \infty$, where $h_{\infty} < 1$.
We note that any constant value of $h_{\infty}$ is allowed
thanks to the invariance of the linear advection--diffusion
equation (\ref{pde}) with respect to the shift and scaling transformations. Indeed, if $h(x,t)$
solves the boundary--value problem (\ref{pde})--(\ref{bc-pde}) such that
$h \to 0$ as $x \to \infty$, then $\tilde{h}(\tilde{x},\tilde{t})$
given by
\begin{equation}
\tilde{h}(\tilde{x},\tilde{t}) = h_{\infty} + (1 - h_{\infty}) h(x,t), \quad
\tilde{x} = (1-h_{\infty})^{1/3} x, \quad \tilde{t} = (1 - h_{\infty})^{4/3} t,
\end{equation}
for any $h_{\infty} < 1$, solves the same advection--diffusion equation (\ref{pde})
with the same boundary conditions (\ref{bc-pde}) but with the variable velocity
$\tilde{V}(\tilde{t}) = \frac{V(t)}{1-h_{\infty}}$ and with the asymptotic behavior
$h \to h_{\infty}$ as $x \to \infty$.

With three boundary conditions at $x = 0$ and the decay conditions for $h$ as $x \to \infty$,
the initial-value problem for equation (\ref{pde}) is over-determined and the third
(over-determining) boundary condition at $x = 0$ is used to find the dependence of $V$ on $t$. Local existence of
solutions to the boundary--value problem (\ref{pde})--(\ref{bc-pde})
was proved by Pelinovsky {\em et al.} \cite{Pelinovsky} using Laplace transform in $x$
and the fractional power series expansion in powers of $t^{1/4}$.

We shall consider the time evolution of the boundary--value problem
(\ref{pde})--(\ref{bc-pde}) starting with the initial data
$h|_{t = 0} = h_0(x)$ for a suitable function $h_0$. For physically relevant solutions,
we assume that the profile $h_0(x)$ decays monotonically to a constant value as $x \to \infty$
and that $0$ is a non-degenerate maximum of $h_0$ such that
$h_0(0) = 1$, $h_0'(0) = 0$, and $h_0''(0) < 0$. The solution $h(x,t)$
may lose monotonicity in $x$ during the dynamical evolution because of the boundary value
\begin{equation}
\label{beta}
\beta(t) := h_{xx}(0,t)
\end{equation}
crosses zero from the negative side. In this case, we say that
the flow becomes non-physical for further times and the model (\ref{pde})--(\ref{bc-pde})
breaks. Simultaneously, this may mean that the velocity $V(t)$ blows up to infinity,
because for sufficiently strong solutions of the advection--diffusion equation (\ref{pde}),
the velocity $V(t)$ satisfies the dynamical equation
\begin{equation}
\label{contact-equation}
h_{xxxxx}(0,t) = V(t) \beta(t),
\end{equation}
which follows by differentiation of (\ref{pde}) in $x$ and setting $x \to 0$.

Based on numerical computations of the thin-film equations (\ref{model}),
Benilov and Vynnycky \cite{Benilov} claim that for any physically relevant
initial data $h_0$, there is a finite positive time $t_0$ such that $V(t)$ tends to negative infinity and $\beta(t)$
approaches zero as $t \nearrow t_0$, whereas the profile $h(x,t_0)$ remains a smooth and decreasing
function for $x > 0$. Moreover, they claim that $V(t)$ behaves near the blowup time as the logarithmic function of $t$:
\begin{equation}
\label{loglaw}
V(t) \sim C_1 \log(t_0 - t) + C_2, \quad \mbox{\rm as} \quad t \nearrow t_0,
\end{equation}
where $C_1$, $C_2$ are positive constants. The same properties of the blow up of contact lines were observed
in \cite{Benilov} in numerical simulations of the reduced model (\ref{pde})--(\ref{bc-pde}). We point out
that the numerical simulations in \cite{Benilov} are based on COMSOL built-in algorithms.

The goal of this paper is to simulate numerically the behavior of the velocity $V(t)$
near the blow-up time under different physically relevant initial data $h_0(x)$.
Our technique is based on the reformulation of the boundary-value problem
(\ref{pde})--(\ref{bc-pde}), which will be suitable for an application
of the direct finite-difference method. We will approximate the behavior of
the velocity $V(t)$ from the dynamical equation (\ref{contact-equation})
rewritten in finite differences. The numerical computations reported in this paper
were performed by using the MATLAB software package.

As the main outcome, we confirm that all physically relevant initial data
including those with positive initial velocity will result in blow-up of $V(t)$
to negative infinity in a finite time. At the same time, we show that the power function
\begin{equation} \label{pwrlaw}
|V(t)| \sim \frac{c}{(t_0 - t)^p}, \quad \text{as} \quad t  \nearrow t_0,
\end{equation}
with $c > 0$ and $p \approx 0.5$ fits our numerical data better
than the logarithmic function \eqref{loglaw} near the blow-up time $t_0$.
We explain why the behavior $|V(t)| = \mathcal{O}((t_0 - t)^{-1/2})$ as $t \nearrow t_0$
is highly expected for solutions of the boundary--value problem (\ref{pde})--(\ref{bc-pde}).
We believe that the incorrect logarithmic law (\ref{loglaw}) is an artefact of the COMSOL
built-in algorithms used in \cite{Benilov}.

We shall mention two recent relevant works on the same problem. Firstly,
existence of self-similar solutions of the linear advection--diffusion equation (\ref{pde})
was proved by Pelinovsky and Giniyatullin \cite{Giniyatullin}. The self-similar solutions
are given by
\begin{equation}
\label{self-similar}
V(t) = \frac{t_0 V_0}{(t_0-t)^{3/4}}, \quad
h(x,t) = f(\xi), \quad \xi = \frac{x}{(t_0-t)^{1/4}},
\end{equation}
with $t_0 > 0$ and $V_0 > 0$, where $f(\xi)$ is a suitable function.
Although the self-similar solutions (\ref{self-similar}) satisfy
the decay condition at infinity, and the first two boundary conditions
(\ref{bc-pde}), the third boundary condition $h_{xxx} |_{x = 0} = - \frac{1}{2}$
is not satisfied and is replaced with $h_{xxx} |_{x = 0} = \gamma_0 V(t)$ for a fixed
$\gamma_0 < 0$. Consequently, the self-similar solution (\ref{self-similar})
predicts blows up in a finite time with positive $V(t)$ and
positive $\beta(t)$. Although the scaling of the self-similar solution (\ref{self-similar})
is compatible with the scaling transformation (\ref{asymptotic-reduction})
used in the derivation of the linear advection--diffusion equation (\ref{pde}),
it does not satisfy the physical requirements of the Couette flow on Figure \ref{fig-scheme}.

Secondly, Chugunova {\rm et al.} \cite{Chugunova} constructed
steady state solutions of the boundary--value problem (\ref{pde})--(\ref{bc-pde})
and showed numerically that these steady states can serve as attractors
of the bounded dynamical evolution of the model. Both the steady states and
the initial conditions that lead to bounded dynamics of the model are not
physically acceptable as $h_0$ has to be monotonically increasing with
$h \to h_{\infty} > 1$ as $x \to \infty$. Note that both $V$ and $\beta$ are
positive for the steady states of the boundary--value problem (\ref{pde})--(\ref{bc-pde}).

To simulate the boundary--value problem (\ref{pde})--(\ref{bc-pde}), a different numerical method
is proposed in \cite{Chugunova}. This method is still based on finite differences and
MATLAB software package. Because the fourth-order derivative term is approximated implicitly and the
first-order derivative term is approximated explicitly, the system of finite-difference
equations was closed in \cite{Chugunova} without any additional equation on the velocity $V(t)$.
Consequently, $V(t)$ was found from the system of finite-difference equations.

We also mention that both recent works of \cite{Giniyatullin} and \cite{Chugunova} used
a priori energy estimates and found some sufficient conditions, under which the smooth
physically relevant solutions of the boundary--value problem (\ref{pde})--(\ref{bc-pde})
blows up in a finite time. In particular, if $V(t) >-1$, or
$\beta(t) < 0$, or $V(t) \beta(t)^2 < 0$, the smooth solution $h(x,t)$
blows up in a finite time. However, these sufficient conditions do not eliminate
existence of smooth physically relevant solutions, for which $V(t)$ oscillates and
$\beta(t)$ decays to zero as $t \to \infty$.

The remainder of our paper is organized as follows. Section 2 outlines the numerical
method for approximations of the boundary--value problem (\ref{pde})--(\ref{bc-pde}).
Section 3 presents the numerical simulations of the boundary--value problem truncated
on the finite interval $[0,L]$ for sufficiently large $L$. Section 4 inspects different
blow-up rates of the singular behavior of the velocity $V(t)$ near the blow-up time.
Section 5 summarizes our findings.

\section{Numerical method}

In what follows, we set $u := h_x$, and reformulate the boundary-value problem \eqref{pde}--\eqref{bc-pde}
in the equivalent form. Differentiating equation \eqref{pde} with respect to $x$, we obtain
\begin{equation} \label{pde2}
\frac{\partial u}{\partial t} + \frac{\partial^4 u}{\partial x^4} = V(t) \frac{\partial u}{\partial x}, \quad x > 0, \quad t > 0.
\end{equation}
We also rewrite boundary conditions in \eqref{bc-pde} as follows:
\begin{equation} \label{bc2}
u|_{x=0} = 0, \quad u_{xx} |_{x=0} = -\frac{1}{2}, \quad u_{xxx} |_{x=0} = 0, \quad t \geq 0.
\end{equation}
Here the third boundary condition $u_{xxx}|_{x=0} = h_{xxxx}|_{x=0} = 0$ follows from
applying the boundary conditions $h|_{x=0} = 1$ and $h_x|_{x=0} = 0$ to the fourth-order equation \eqref{pde}
as $x \rightarrow 0$. After the reformulation, the dynamical equation
(\ref{contact-equation}) can be recovered by taking the limit $x \rightarrow 0$ in \eqref{pde2}:
\begin{equation} \label{velocity1}
u_{xxxx}(0,t) = V(t) u_x(0,t), \quad t \geq 0,
\end{equation}
provided that the solution $u$ remains smooth at the boundary $x = 0$.

A suitable two-parameter initial condition $u|_{t=0} = u_0(x)$ for the boundary--value problem
\eqref{pde2}--\eqref{bc2} can be chosen in the form
\begin{equation} \label{ic}
u_0(x) = -\frac{1}{4}e^{-ax}x[4 + (4a+1)x + a(2a+1)x^2 + bx^3],
\end{equation}
where parameters $a > 0$ and $b \geq 0$ are arbitrary. For simplicity,
the constraint
$$
\beta(t) = h_{xx}|_{x=0} = u_x|_{x=0} < 0
$$
at the initial time $t = 0$ is standardized to $\beta(0) = -1$.
Note that $u_0(x)$ and its derivatives decay to zero exponentially fast as $x \to \infty$,
which still imply that $h_0(x) = 1 + \int_0^x u_0(x') dx'$ decays to some
constant value $h_{\infty}$ as $x \to \infty$. Because $u_0(x) < 0$ for all $x > 0$,
$h_0$ is a monotonically decreasing function and $h_{\infty} < 1$.

\begin{figure}[H]
\includegraphics[height=2in]{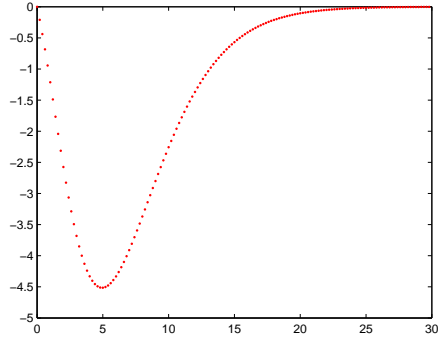}
\caption{The initial function \eqref{ic} with $a = 0.5$ and $b = 0$.}
\label{fig-initial}
\end{figure}

Figure \ref{fig-initial} shows a particular example of the initial function (\ref{ic})
for $a = 0.5$ and $b = 0$.

We approximate solutions of the boundary--value problem (\ref{pde2})--(\ref{bc2}) with the second-order
central-difference method. Consider a set of $N+2$ equally spaced ordered grid points
$\lbrace x_n \rbrace _{n=0}^{N+1}$ on the interval $[0, L]$, for sufficiently large $L$
so that $u|_{x=L}$ and $u_{xx} |_{x = L}$ are approximately zero. For any fixed $t > 0$, let $u_n(t)$
denote the numerical approximation of $u(x,t)$ at $x = x_n$, and let $\Delta x$ denote
the equal step size between adjacent grid points.

By applying the second-order central-difference formulas to partial derivatives in
the fourth-order equation \eqref{pde2} at each $x = x_n$, we obtain the differential equations:
\begin{equation} \label{fd1}
\frac{du_n}{dt} = V(t) \frac{u_{n+1} - u_{n-1}}{2(\Delta x)} - \frac{u_{n+2} - 4u_{n+1} + 6u_n - 4u_{n-1} + u_{n-2}}{(\Delta x)^4},
\end{equation}
which are accurate up to the $\mathcal{O}(\Delta x^2)$ truncation error.
Since $u_0 = u(0, t) = 0$ and $u_{N+1} = u(L, t) = 0$ for all $t \geq 0$, the above formula needs
only to be applied to $N$ interior points $\lbrace x_n \rbrace_{n=1}^N$ with the necessity to
approximate $u_{-1}$ for the grid point $x_1$ and $u_{N+2}$ for the grid point $x_N$.
The value of $u_{-1}$ can be found from the boundary condition:
\begin{equation*}
u_{xx} |_{x=0} = -\frac{1}{2} \quad \Rightarrow \quad
\frac{u_{-1} - 2u_0 + u_1}{(\Delta x)^2} = -\frac{1}{2} \quad \Rightarrow \quad u_{-1} = -u_1 - \frac{1}{2}(\Delta x)^2,
\end{equation*}
and $u_{N+2}$ can be found from the decay condition:
\begin{equation*}
u_{xx} |_{x = L} = 0 \quad \Rightarrow \quad \frac{u_N - 2u_{N+1} + u_{N+2}}{(\Delta x)^2} = 0
\quad \Rightarrow \quad u_{N+2} = -u_N,
\end{equation*}
which are again accurate up to the $\mathcal{O}(\Delta x^2)$ truncation error.
It remains to define $V(t)$ from the third boundary condition $u_{xxx} |_{x=0} = 0$.

The velocity $V(t)$ can be expressed by applying the central--difference approximation
to the dynamical equation \eqref{velocity1}:
\begin{equation} \label{velocity2}
V(t) = \frac{u_{xxxx}|_{x=0}}{u_x|_{x=0}} \quad \Rightarrow \quad
V(t) = \frac{2(u_2 - 4u_1 + 6u_0 - 4u_{-1} + u_{-2})}{(\Delta x)^3(u_1 - u_{-1})},
\end{equation}
where $u_{-2}$ can be found from the third boundary condition in \eqref{bc2}:
\begin{equation*}
u_{xxx} |_{x=0} = 0 \quad \Rightarrow \quad
\frac{u_2 - 2u_1 + 2u_{-1} - u_{-2}}{(\Delta x)^3} = 0 \quad \Rightarrow \quad u_{-2} = u_2 - 4u_1 - (\Delta x)^2.
\end{equation*}

Writing the system of differential equations \eqref{fd1} in the matrix form
$$
\frac{d\pmb{u}}{dt} = \pmb{A}\pmb{u} + \pmb{b},
$$
we use Heun's method to evaluate solutions of the system of differential equations.
Let $\pmb{u}_k$ denote the numerical approximation of $\pmb{u}(t)$ at $t = t_k$
and let $\Delta t$ denote the time step size (not necessarily constant).
By Heun's method, we obtain the iterative rule
\begin{equation} \label{heun1}
\pmb{u}_{k+1} = \pmb{u}_k + \frac{\Delta t}{2}[(\pmb{A}_k \pmb{u}_k + \pmb{b}) + (\pmb{A}_{k+1} \pmb{u}_{k+1} + \pmb{b})],
\end{equation}
where the initial vector $\pmb{u}_0$ is obtained from the initial condition \eqref{ic}. Note that the
coefficient matrix $\pmb{A}$ depend on $t$ since it is defined by the variable velocity $V(t)$.
Nevertheless, ${\bf b}$ is constant in $t$.
The global error of Heun's method is $\mathcal{O}(\Delta t^2)$, so the global truncation error
for the numerical approximation is $\mathcal{O}(\Delta x^2 + \Delta t^2)$.

The explicit version of Heun's method is stable only when
$$
\Delta t \leq \displaystyle\frac 1 8 (\Delta x)^4.
$$
Therefore, in practice we shall use the implicit Heun's method (which is stable for all $\Delta t > 0$),
by solving the system of linear equations
\begin{equation} \label{heun2}
(\pmb{I} - \frac{\Delta t}{2}\pmb{A}_{k+1})\pmb{u}_{k+1} = (\pmb{I} + \frac{\Delta t}{2}\pmb{A}_{k})\pmb{u}_k + \Delta t \pmb{b},
\end{equation}
where $\pmb{I}$ is the identity matrix. However, because the coefficient matrix $\pmb{A}_{k+1}$ on the left-hand
side contains an unknown value of $V(t_{k+1})$, a prediction-correction method is necessary for solving this system
of equations as follows. First, $\pmb{A}_{k+1}$ is approximated using $\pmb{A}_k$ to predict the value of
$\pmb{u}_{k+1}^*$, which is then used to predict the value of $V(t_{k+1}^*)$ using equation
\eqref{velocity2}. Second, $\pmb{A}_{k+1}$ is updated from the prediction $V(t_{k+1}^*)$ to obtain
the corrected values of $\pmb{u}_{k+1}$ and $V(t_{k+1})$. Since the implicit method is used in
both prediction and correction steps, the unconditional stability is preserved.

\section{Blow-up of the velocity $V$ of contact lines}

We use the finite-difference method to compute approximation of the boundary--value problem
(\ref{pde2})--(\ref{bc2}), after truncation on the finite interval $[0,L]$ with sufficiently
large $L$. Since the time evolution features blow-up in a finite time,
an adaptive method is used to adjust the time step $\Delta t$ after each iteration to maintain
the local truncation error of the numerical method at a certain tolerance level.

\begin{figure}[H]
\begin{center}
\includegraphics[height = 2in]{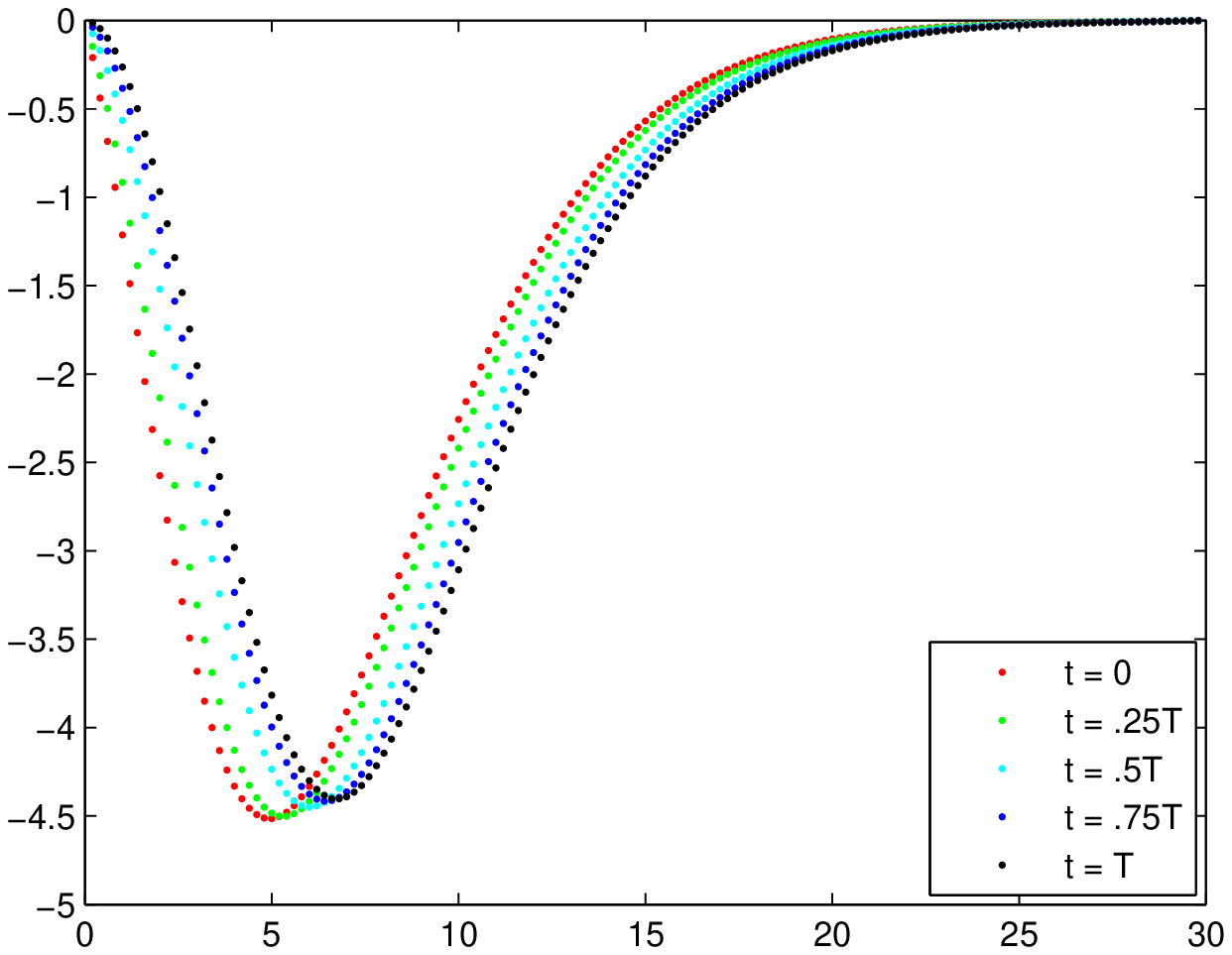}
\includegraphics[height = 2in]{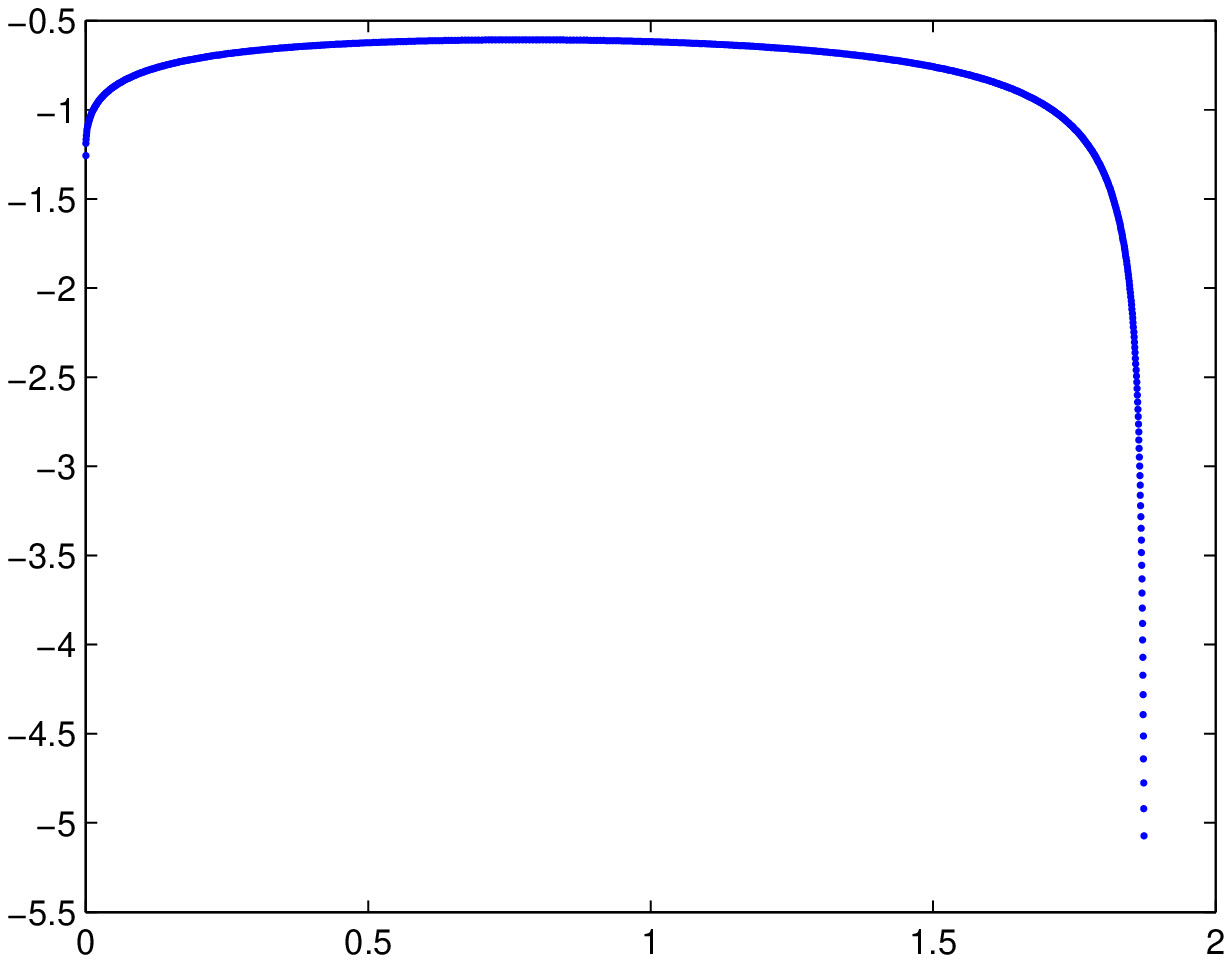} \\
\includegraphics[height = 2 in]{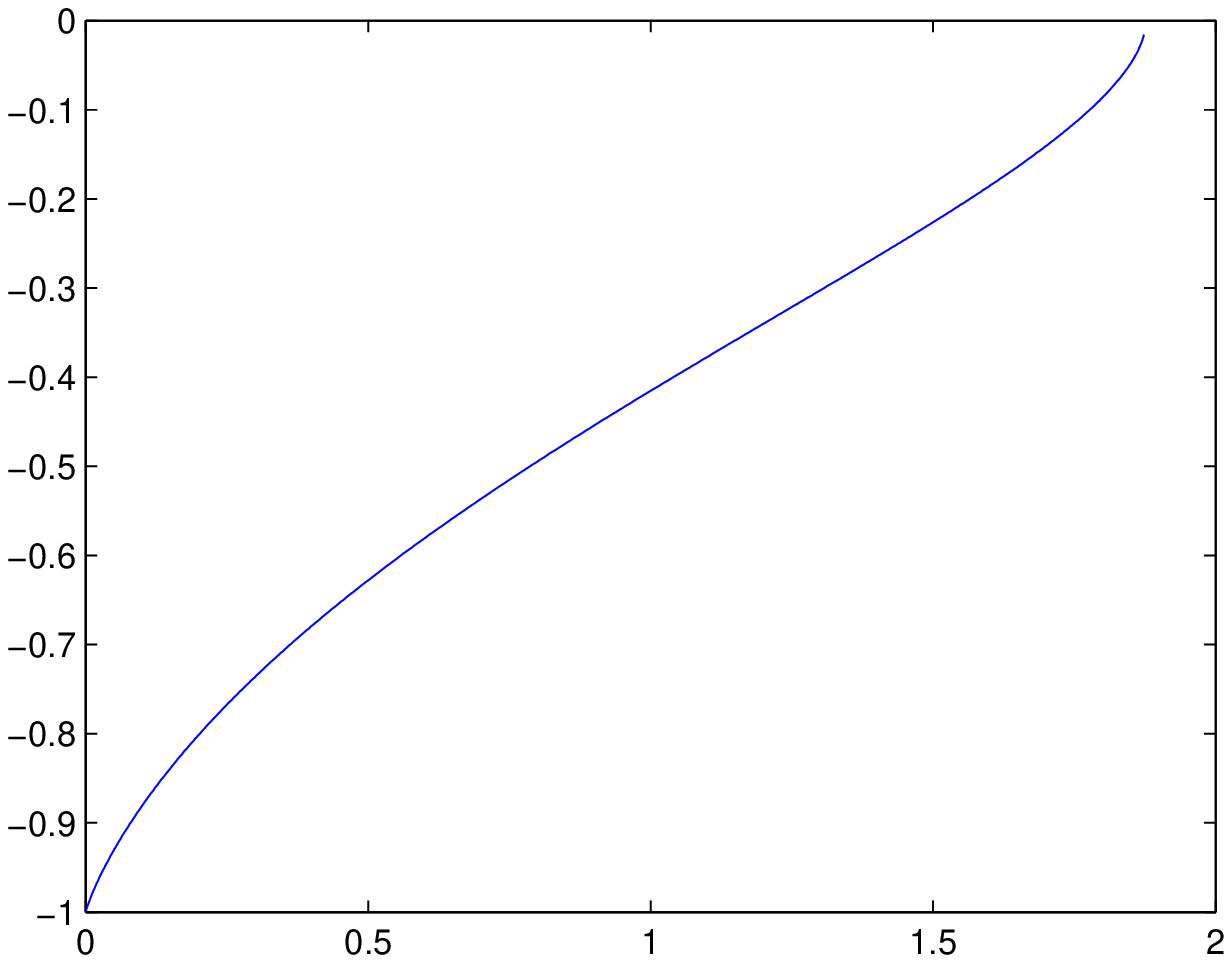}
\includegraphics[height = 2 in]{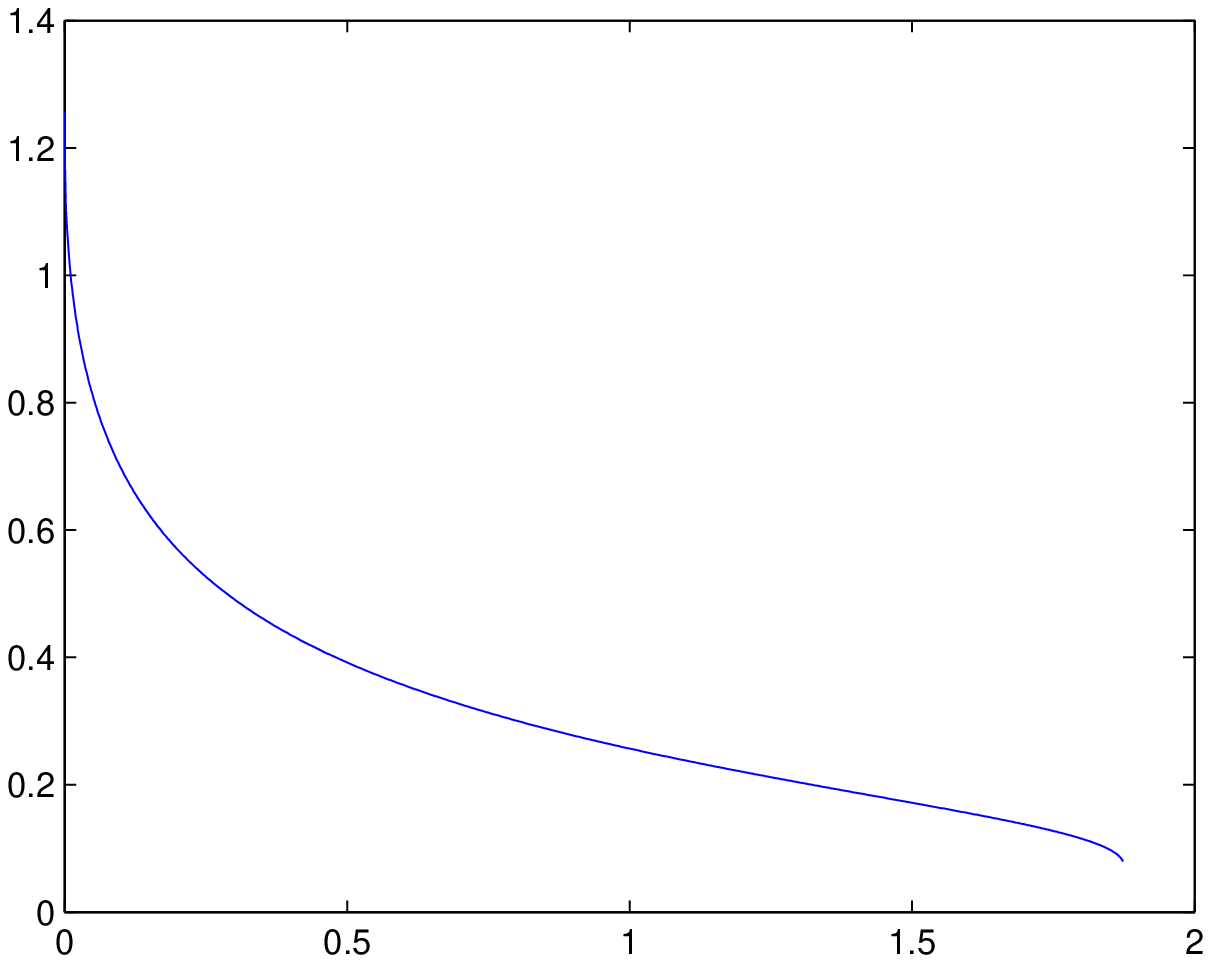}
\caption{Numerical solution of the boundary-value problem (\ref{pde2})--(\ref{bc2}),
where $u_0$ is given by (\ref{ic}) with $a = 0.5$ and $b = 0$.
Top left: The profile of $u$ versus $x$ at different time $t$.
Top right: Velocity of the contact line $V$ versus $t$. Bottom left: The
boundary value $\beta = u_x |_{x = 0}$ versus $t$. Bottom right: The
boundary value $u_{xxxx} |_{x = 0}$ versus $t$.}
\label{fig-num-1}
\end{center}
\end{figure}

Figure \ref{fig-num-1} shows the numerical approximation of the boundary--value problem
(\ref{pde2})--(\ref{bc2}) for the initial function \eqref{ic} with $a = 0.5$ and $b = 0$
(the one shown on Figure \ref{fig-initial}). The initial velocity is determined from
this initial condition by equation (\ref{velocity1})
as $V(0) = -1.25$. The top left panel of the figure shows the profile of $u(x,t)$ versus $x$
at different values of $t$ until the terminal time $T$ of our computations.
The top right panel of the figure shows the change of the velocity $V(t)$ in time $t$
computed dynamically from equation (\ref{velocity2}). The bottom left panel shows the
boundary value $\beta(t) = u_x |_{x = 0}$ versus $t$ and the bottom right panel shows the
boundary value $u_{xxxx} |_{x = 0}$ versus $t$.

It is clear from the top panels that the velocity $V$ diverges towards $-\infty$ at
$t \approx 1.9$, whereas the solution $u(x,t)$ remains regular near the blow-up time.
Recall that the velocity $V(t)$ is determined from equation \eqref{velocity1}
by the quotient of $u_{xxxx}(0,t)$ and $\beta(t) = u_x(0,t)$,
where $\beta(t)$ must be strictly negative for all $t > 0$ for physically acceptable solutions.
We can see from the bottom panels that the value of $\beta$ is about to cross zero from the negative side
at the blow-up time, whereas $u_{xxxx}(0, t)$ is also approaching zero but at a much slower rate than $\beta(t)$.
This also indicates that $V(t)$ is approaching negative infinity at the blow-up time.

To measure the error of numerical computations, we shall derive dynamical constraints
on the time evolution of a smooth solution of the boundary--value problem (\ref{pde2})--(\ref{bc2}).
Differentiating equation \eqref{pde2} with respect to $x$ once and twice and taking the limit
$x \rightarrow 0$, we obtain
\begin{equation}
\label{r5}
\frac{d\beta}{dt} + \frac{\partial^5 u}{\partial x^5}\biggr|_{x=0} = -\frac{1}{2}V(t)
\end{equation}
and
\begin{equation}
\label{r6}
\frac{\partial^6 u}{\partial x^6}\biggr|_{x=0} = 0.
\end{equation}
Using equation \eqref{r6}, we determine $u_{-3}$ at $x = x_{-3}$ from the central--difference approximation:
\begin{eqnarray*}
& \phantom{t} & \frac{u_3 - 6u_2 + 15u_1 - 20u_0 + 15u_{-1} - 6u_{-2} + u_{-3}}{(\Delta x)^6} = 0\\
& \phantom{t} & \phantom{texttexttexttexttexttext}
\Rightarrow \quad u_{-3} = -u_3 + 12u_2 - 24u_1 + \frac{3}{2}(\Delta x)^2.
\end{eqnarray*}
Then, the value of $d\beta/dt$ is approximated from equations \eqref{velocity2} and \eqref{r5}:
\begin{eqnarray}
\beta'(t) = -\frac{u_3 - 4u_2 + 5u_1 - 5u_{-1} + 4u_{-2} - u_{-3}}{2(\Delta x)^5}
- \frac{u_2 - 4u_1 + 6u_0 - 4u_{-1} + u_{-2}}{(\Delta x)^3(u_1 - u_{-1})}. \label{dbeta}
\end{eqnarray}
Comparing the value of $\beta'(t)$ determined from equation \eqref{dbeta} with the central-difference
approximation for the numerical derivative
\begin{equation} \label{dbetanum}
\beta'(t_l) = \frac{\beta(t_{k+1}) - \beta(t_{k-1})}{t_{k+1} - t_{k-1}},
\end{equation}
we can estimate the numerical error of the solution at the boundary $x = 0$.
\begin{figure}[H]
\begin{center}
\includegraphics[height = 2in]{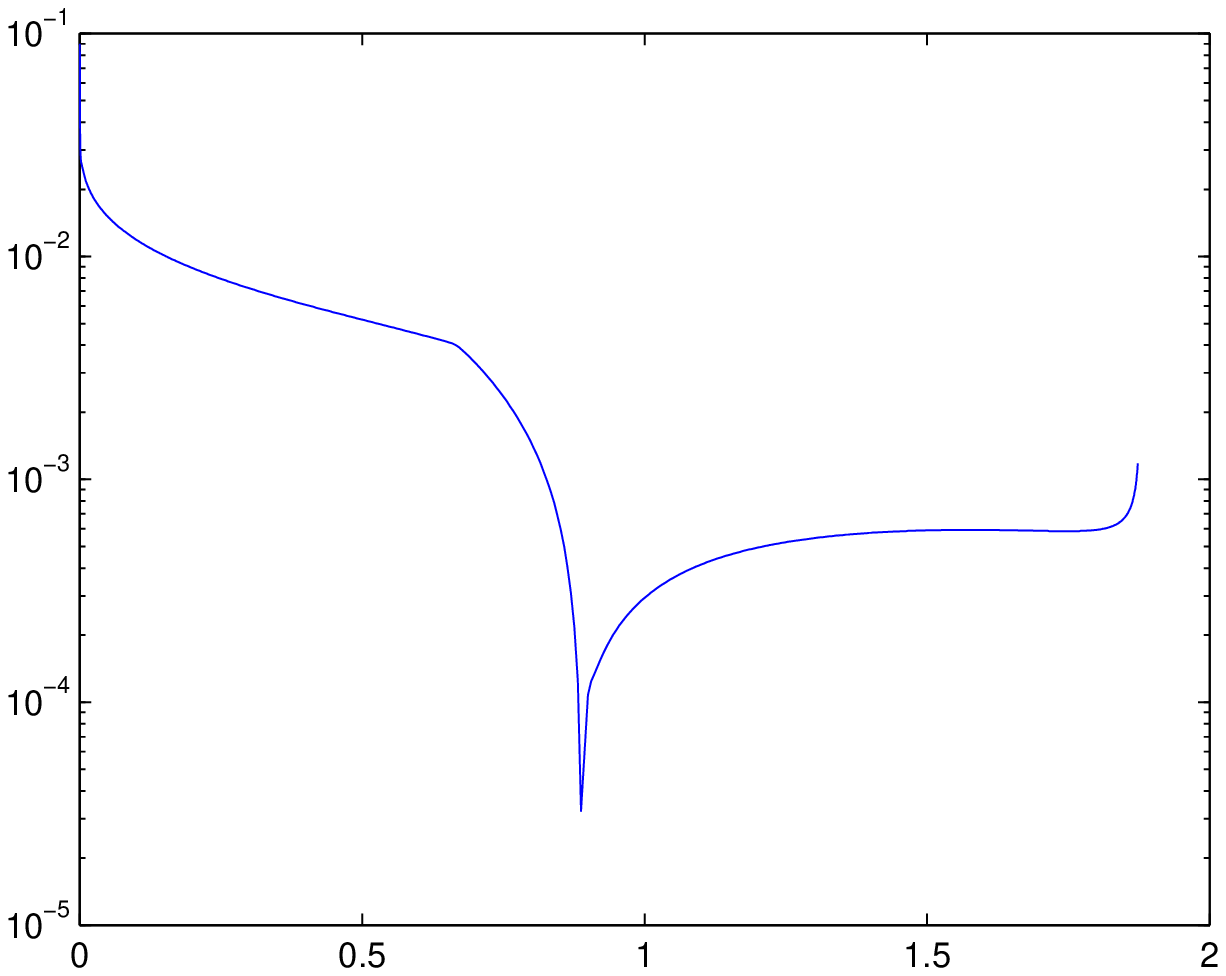}
\includegraphics[height = 2in]{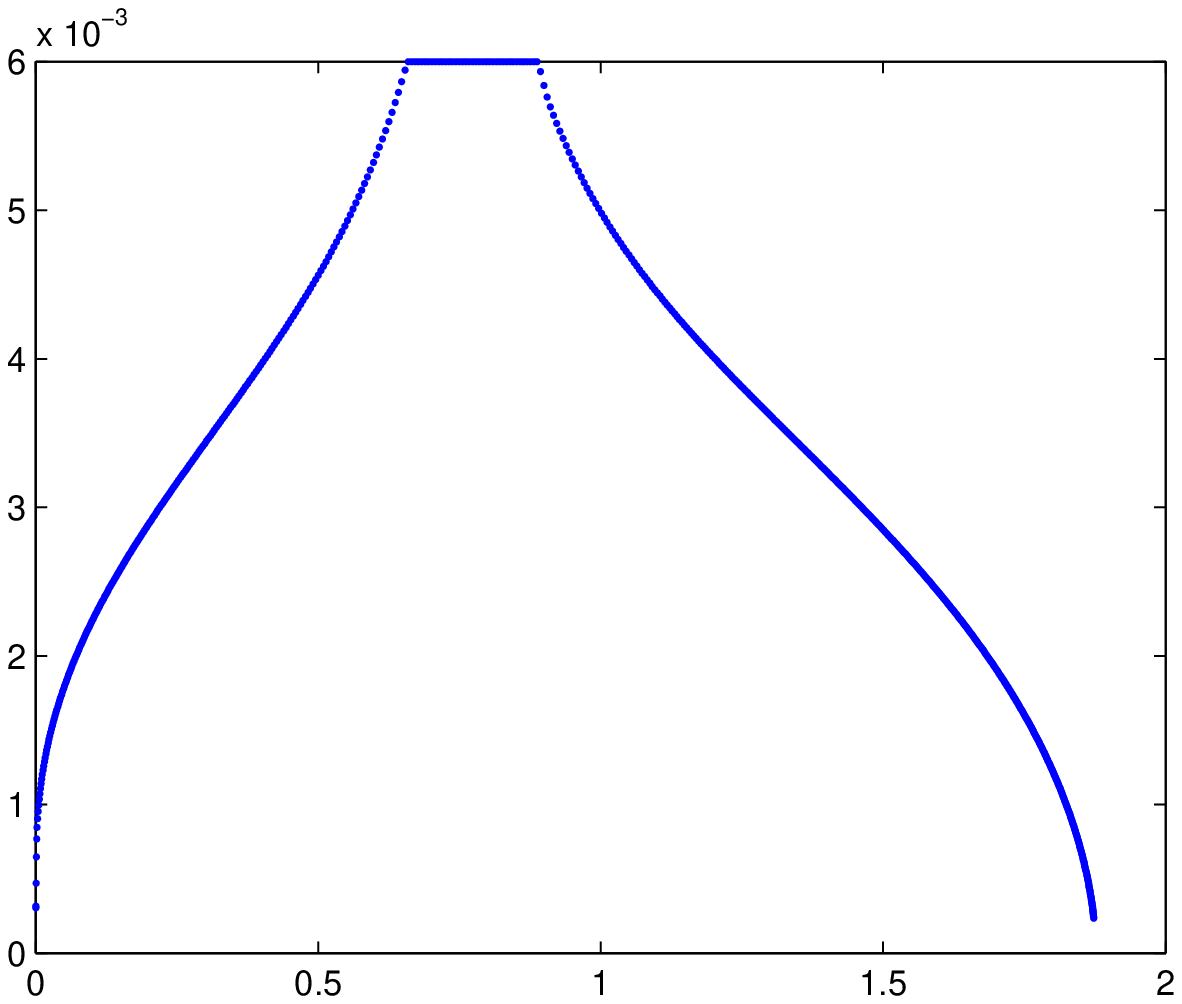}
\caption{Left: Error of $\beta'(t)$ versus $t$.
Right: Time step size $\Delta t$ versus $t$.}
\label{fig-num-2}
\end{center}
\end{figure}

Figure \ref{fig-num-2} (left) compares the value of $\beta'(t)$ between equations
\eqref{dbeta} and \eqref{dbetanum}. The error remains small, therefore, the
assumption that the solution is smooth (or at least $C^6$) at the boundary $x = 0$
is valid up to numerical accuracy.

Figure \ref{fig-num-2} (right) shows the time step size $\Delta t$ adjusted
to preserve the same level of the local error of $10^{-5}$. We set $\Delta t = 0.006$
if the error estimation procedure yields larger values of $\Delta t$.
This truncation is needed because the error drops significantly near $t \approx 0.8$,
and the error estimation procedure would otherwise produce large values of $\Delta t$.

\begin{figure}[H]
\begin{center}
\includegraphics[height = 2in]{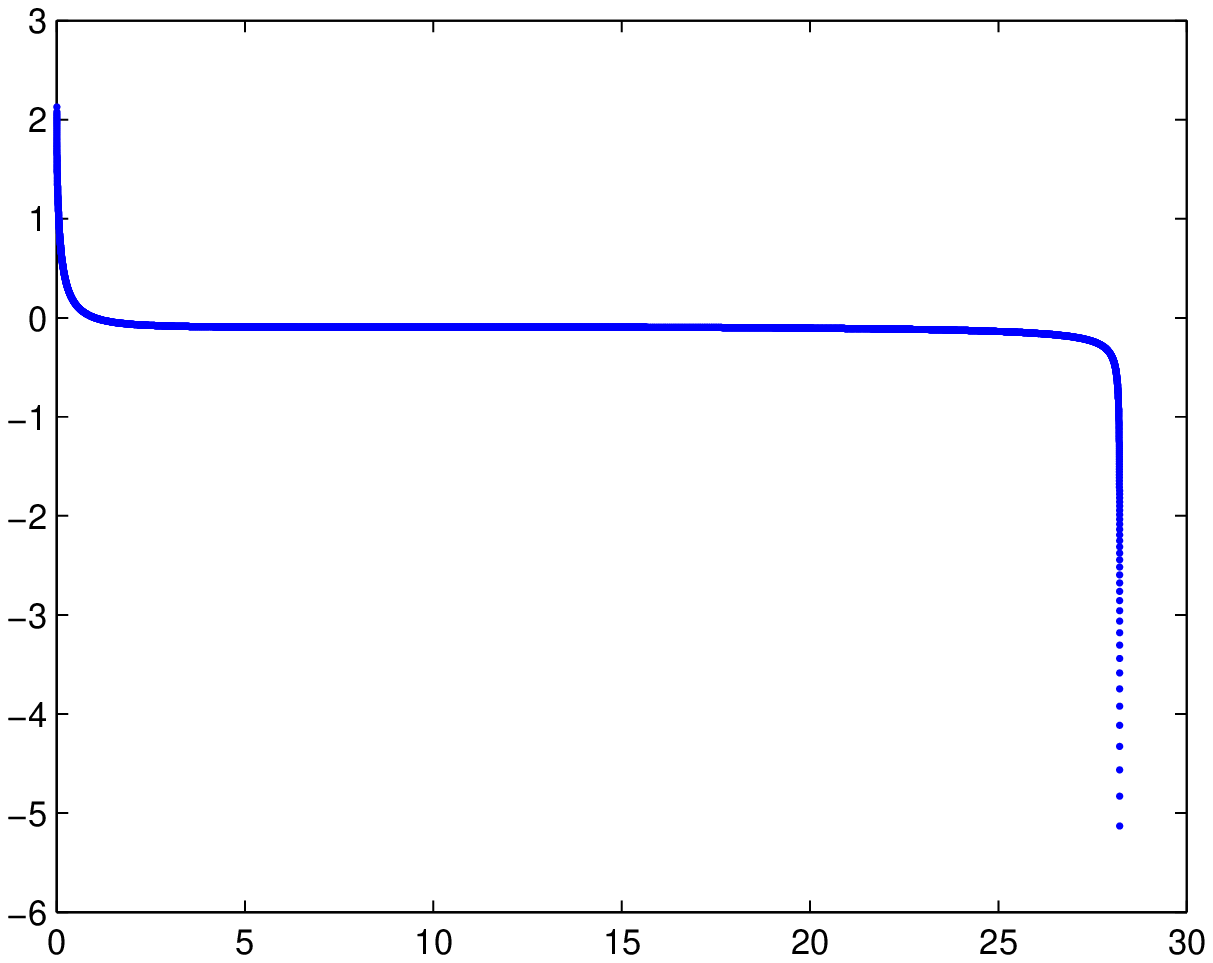}
\includegraphics[height = 2in]{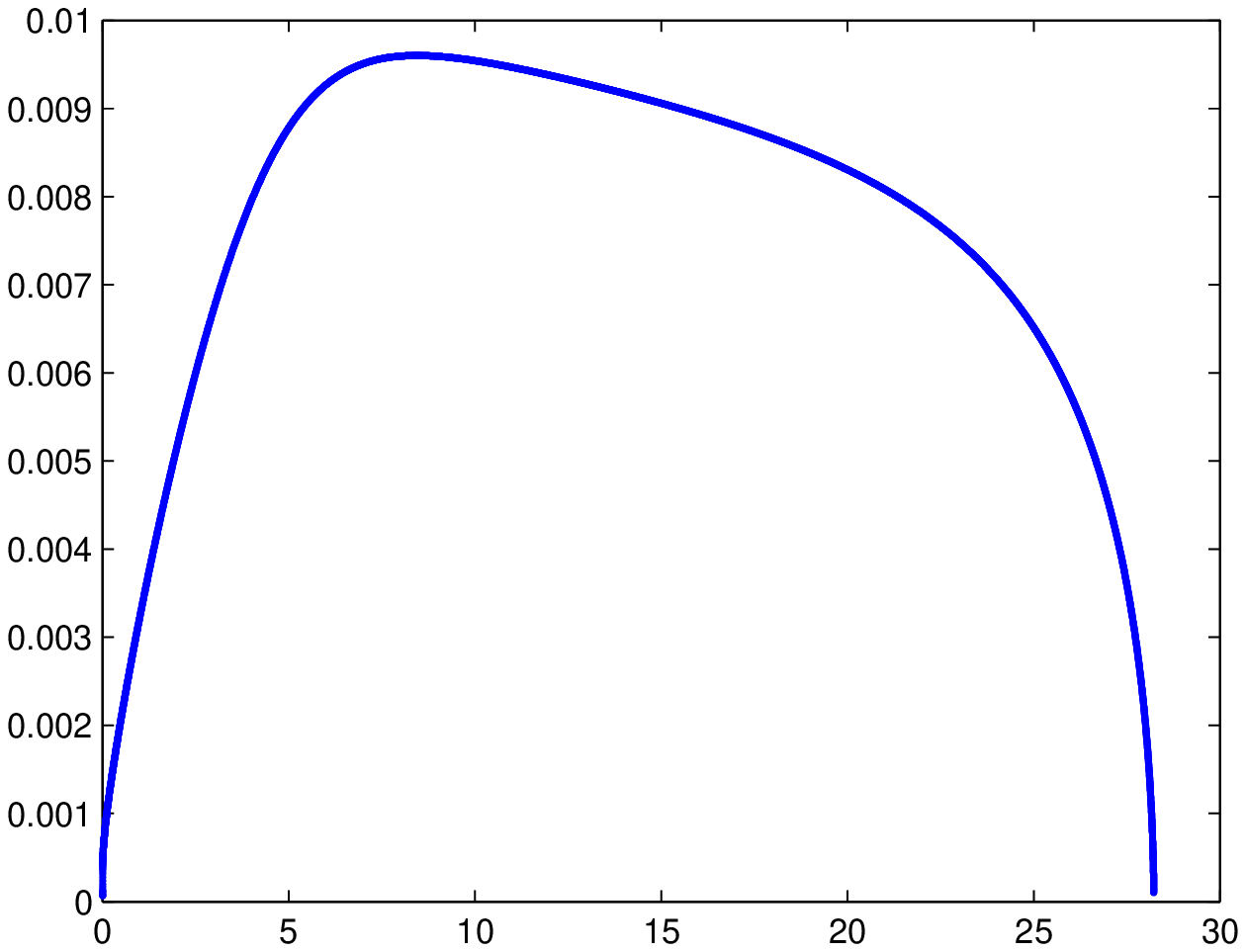}
\caption{Numerical solution of the boundary-value problem (\ref{pde2})--(\ref{bc2}), where
$u_0$ is given by (\ref{ic})
with $a = 0.5$ and $b = 0.6$. Left: Velocity of the contact line $V$ versus $t$.
Right: Time step size $\Delta t$ versus $t$. }
\label{fig-num-3}
\end{center}
\end{figure}

We have performed computations with other initial conditions from the two-parameter family
of functions in (\ref{ic}). Figure \ref{fig-num-3} (left) shows the dynamical evolution
of the velocity $V$ starting with a positive velocity $V(0) = 2.35$, which is
determined from the initial function (\ref{ic}) with $a = 0.5$ and $b = 0.6$.
Although the terminal time $T \approx 28$ is much larger compared with the case of the
negative initial velocity on Figure \ref{fig-num-1}, a blow-up is still detected
from this initial condition. The solution $u(x,t)$ looks similar to the solution shown in
Figure \ref{fig-num-1} (top left) and hence is not shown.

Figure \ref{fig-num-3} (right) shows
the adjusted time step size. We note that the time step size is small at the initial time because
the smooth solution appears from the initial condition, which does not satisfy infinitely many constraints
of the boundary--value problem (\ref{pde2})--(\ref{bc2}). It is also small near the terminal time because
of the blow-up of the smooth solution. But $\Delta t$ is not too small at intermediate values of $t$,
when the solution is at a slowly varying phase. During this slowly varying phase, $V(t)$ is nearly constant but
$\beta(t)$ changes nearly linearly in time (similarly to Figure \ref{fig-num-1} (bottom left)
and hence is not shown).

\begin{figure}[H]
\begin{floatrow}
\includegraphics[height = 2.2 in]{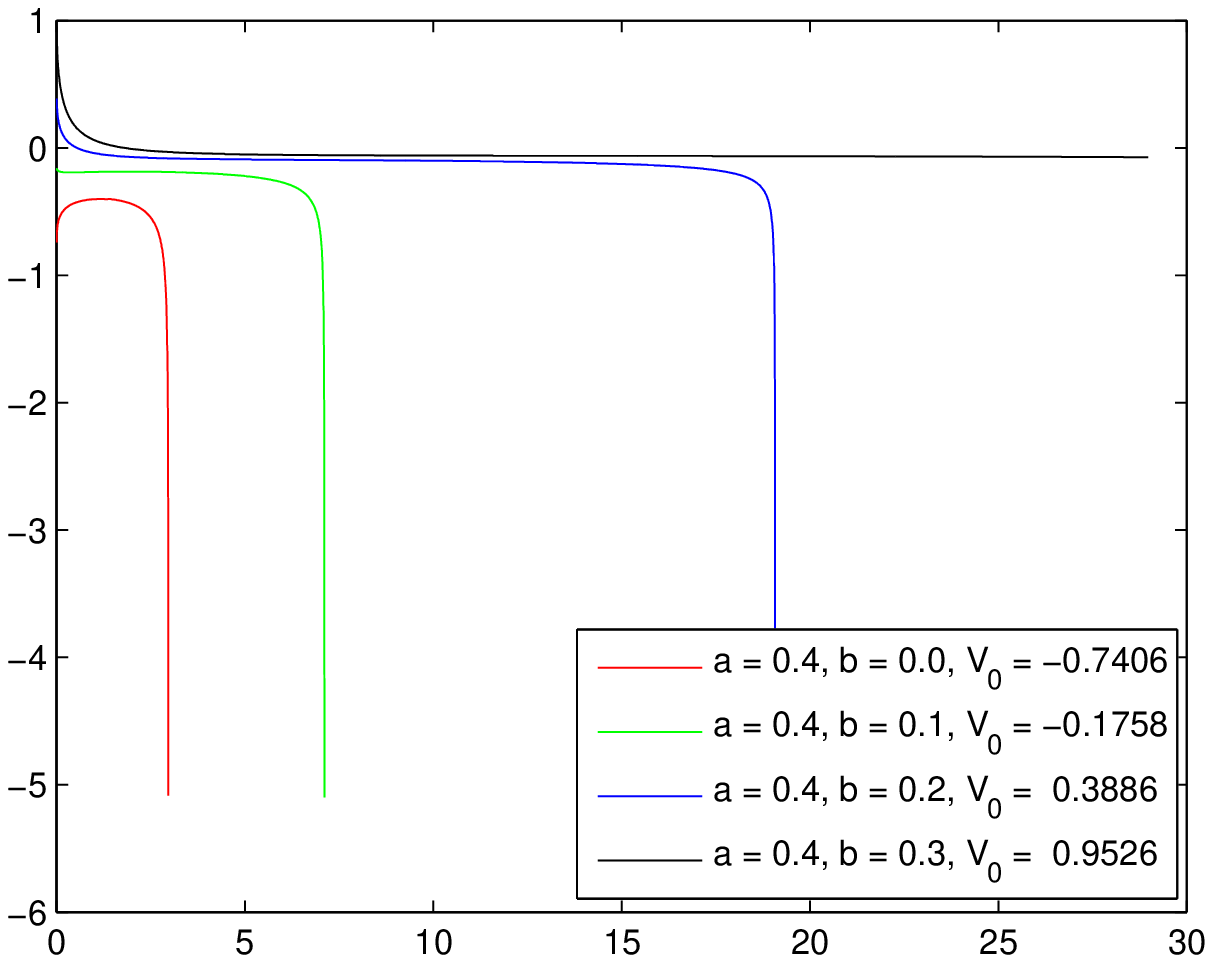}
\includegraphics[height = 2.2 in]{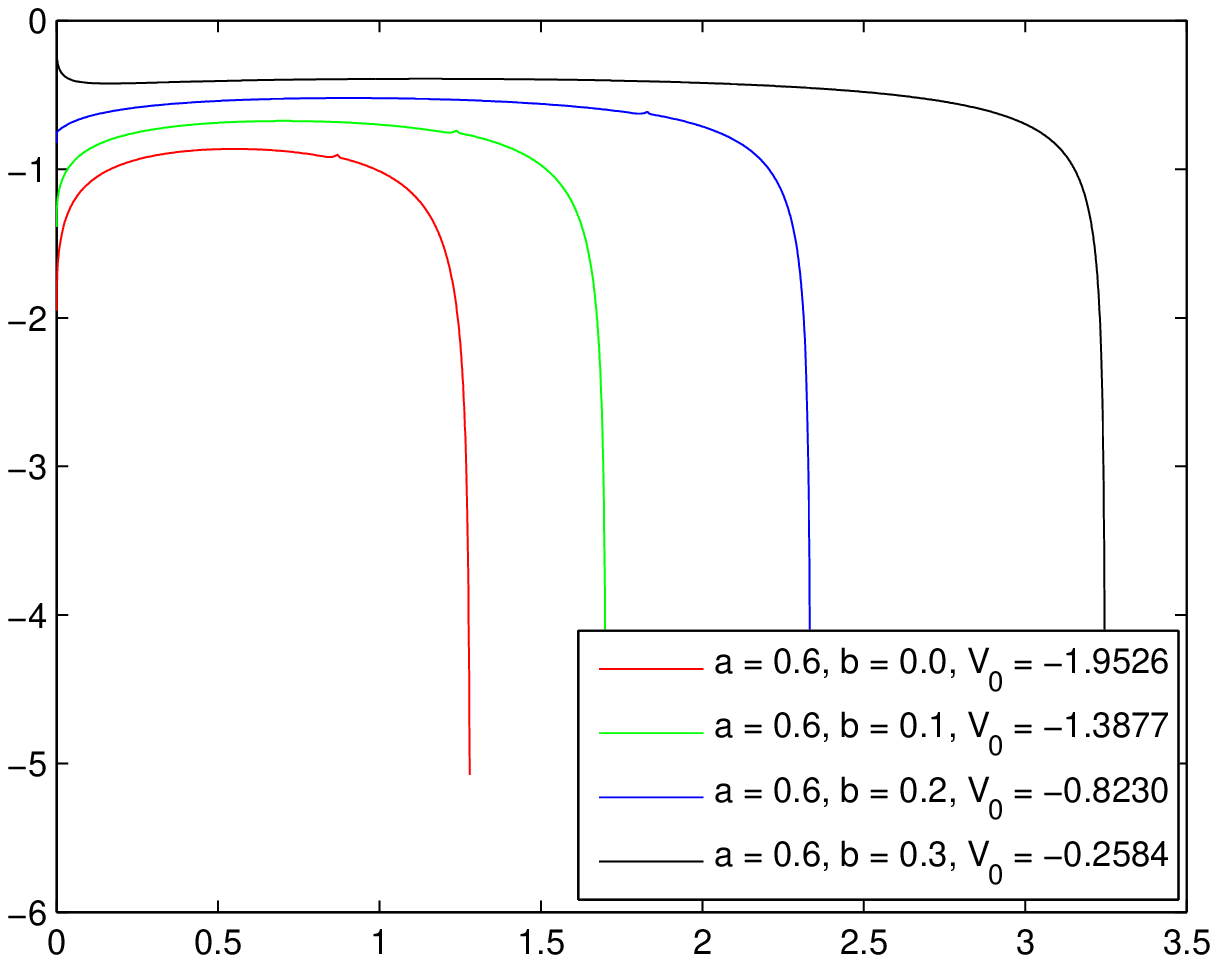}
\end{floatrow}
\caption{Behavior of the velocity $V$ versus $t$ for different initial conditions
given by the two-parameter function (\ref{ic}).}
\label{fig-num-4}
\end{figure}

Figure \ref{fig-num-4} illustrates the dynamical evolution of the velocities $V(t)$ under
different initial conditions given by the two-parameter function (\ref{ic}).
From these plots, together with the previous examples, it is clear that the blow-up time
depends on the initial velocity $V(0)$ and a large positive initial velocity leads to
a much longer slowly varying phase before the solution blows up. Nevertheless, the blow-up
in a finite time is unavoidable for any physically acceptable initial conditions.

\section{Rate of blow-up}

In order to determine numerically the blow-up time $t_0$ and the rate of blow-up of the velocity $V$,
we will fit the numerical data near the terminal time $T$ of our computations
with either the logarithmic function \eqref{loglaw} or the power function \eqref{pwrlaw}.

For the logarithmic function \eqref{loglaw}, we first differentiate both sides
of the expression with respect to $t$ and take the inverse:
\begin{equation}
\label{regression-1}
\frac{dV}{dt} = -\frac{C_1}{t_0 - t} \quad \Rightarrow \quad
\biggr(\frac{dV}{dt}\biggr)^{-1} = \frac{t}{C_1} - \frac{t_0}{C_1}.
\end{equation}
Then the constants $C_1$ and $t_0$ can be determined from a linear regression applied
to equation (\ref{regression-1}). 
We will skip the numerical procedure for determining the values of $C_2$ since it does not affect the 
blow-up behavior of $V$.

For the power function \eqref{pwrlaw}, we can take the logarithm of both sides of the expression
$$
\log(-V(t)) = \log c - p\log(t_0 - t)
$$
and then differentiate the above expression:
\begin{equation}
\label{regression-2}
\frac{1}{V(t)}\frac{dV}{dt} = \frac{p}{t_0 - t} \quad \Rightarrow \quad
V(t)\biggr(\frac{dV}{dt}\biggr)^{-1} = \frac{t_0}{p} - \frac{t}{p}.
\end{equation}
The constants $p$ and $t_0$ can now be determined from a linear regression applied
to equation (\ref{regression-2}).

In practice, we found that the blow-up rate $p$ in the power law
or the coefficient $C_1$ in the logarithmic law vary with different time windows (i.e.
the range of $t$ which is used to fit the data). The following output gives
a comparison of numerical data under different time windows and different
tolerance levels, using the initial condition (\ref{ic}) with $a = 0.5$ and $b = 0$.
Here \emph{starting time} means the time $t$ at which we start to fit the data, and
\emph{Error} is the mean squared error (MSE) defined by
$$
\text{MSE} := \frac{1}{n-3}\sum(V_{\text{obs}} - V_{\text{fit}})^2,
$$
where $n$ is the total number of data points used in the regression.

\begin{verbatim}
Initial condition: a = 0.5, b = 0; initial velocity: V(0) = -1.2500

Tolerance level: 0.0001, number of iterations: 330, terminal time = 1.8729
  Starting time    Blowup time t0    Blowup rate p or C1      Error
powerlaw:
         1.8176         1.8749            0.3916           0.000017
         1.8356         1.8752            0.3994           0.000003
         1.8550         1.8756            0.4104           0.000000
loglaw:
         1.8176         1.8678            0.5371          23.732740
         1.8356         1.8695            0.6135          33.681247
         1.8550         1.8716            0.7578          68.934686

Tolerance level: 1e-006, number of iterations: 1448, terminal time = 1.8732
  Starting time    Blowup time t0    Blowup rate p or C1      Error
powerlaw:
         1.8172         1.8753            0.3927           0.000033
         1.8360         1.8757            0.4009           0.000006
         1.8547         1.8760            0.4118           0.000000
loglaw:
         1.8172         1.8688            0.5500          25.226547
         1.8360         1.8705            0.6343          33.937325
         1.8547         1.8724            0.7854          58.894321
\end{verbatim}

\begin{figure}[H]
\includegraphics[height = 2.5in]{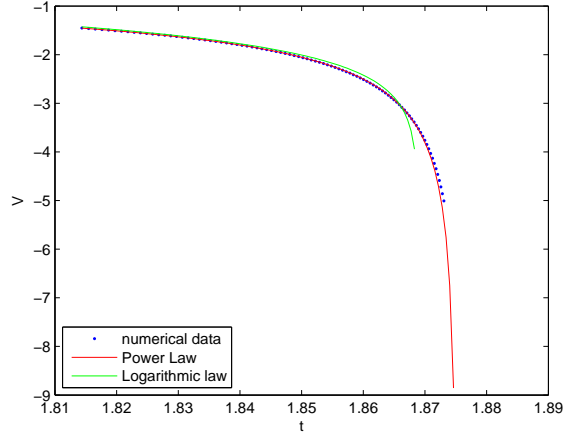}
\caption{Comparison between data fitting with the logarithmic law \eqref{loglaw}
and the power law \eqref{pwrlaw}.}
\label{fg:datafitting}
\end{figure}

The above table shows that the errors from the logarithmic law \eqref{loglaw}
are much larger than the errors from the power law \eqref{pwrlaw} in all cases.
Also, the error of the power law reduces as the time window moves closer to
the blow-up time, whereas the error of the logarithmic law increases. Moreover,
the blow-up times $t_0$ determined from the logarithmic law are smaller than
the terminal time $T$ of our computations. Hence, the logarithmic law deviates
from the numerical data near the blow-up time. As we can see in Figure \ref{fg:datafitting},
the power function \eqref{pwrlaw} fits our numerical data much better than
the logarithmic function \eqref{loglaw}.

In order to confirm that the blow-up of the velocity $V$ occurs according to the power
law (\ref{pwrlaw}) compared to the logarithmic law (\ref{loglaw}), we use the scaling
transformations suggested in \cite{Benilov} and replace the time variable $t$ by
the new variable
\begin{equation}
\label{time-variable}
T := \int_0^t (1 + V^{2n}(t')) dt',
\end{equation}
where $n$ is a positive integer. In new time variable with $V(t) \equiv V(T)$, the model (\ref{pde2})
is rewritten in the form
\begin{equation} \label{pde3}
\frac{\partial u}{\partial T} = \frac{1}{1 + V^{2n}} \left( V \frac{\partial u}{\partial x} -
\frac{\partial^4 u}{\partial x^4}\right), \quad x > 0, \quad T > 0,
\end{equation}
whereas the boundary conditions or the numerical method are unaffected. With the power
law (\ref{pwrlaw}) as $t  \nearrow t_0$, the new time variable $T$ in (\ref{time-variable}) approaches
a finite limit if $2np < 1$ and becomes infinite if $2np \geq 1$. With the logarithmic law
(\ref{loglaw}), the new time variable $T$ would always approach a finite limit for any integer $n$.

\begin{figure}[H]
\begin{floatrow}
\includegraphics[height = 2.2 in]{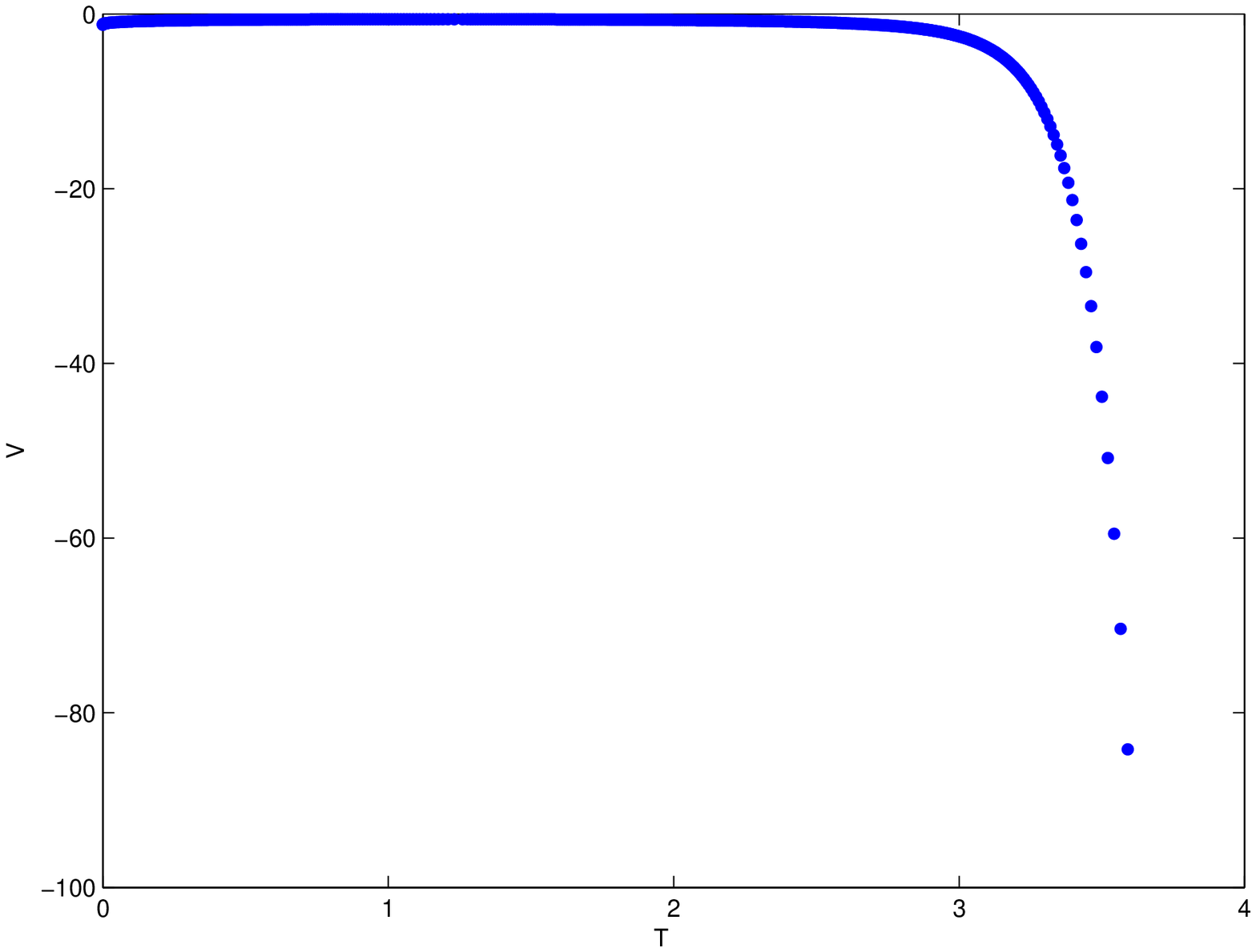}
\includegraphics[height = 2.2 in]{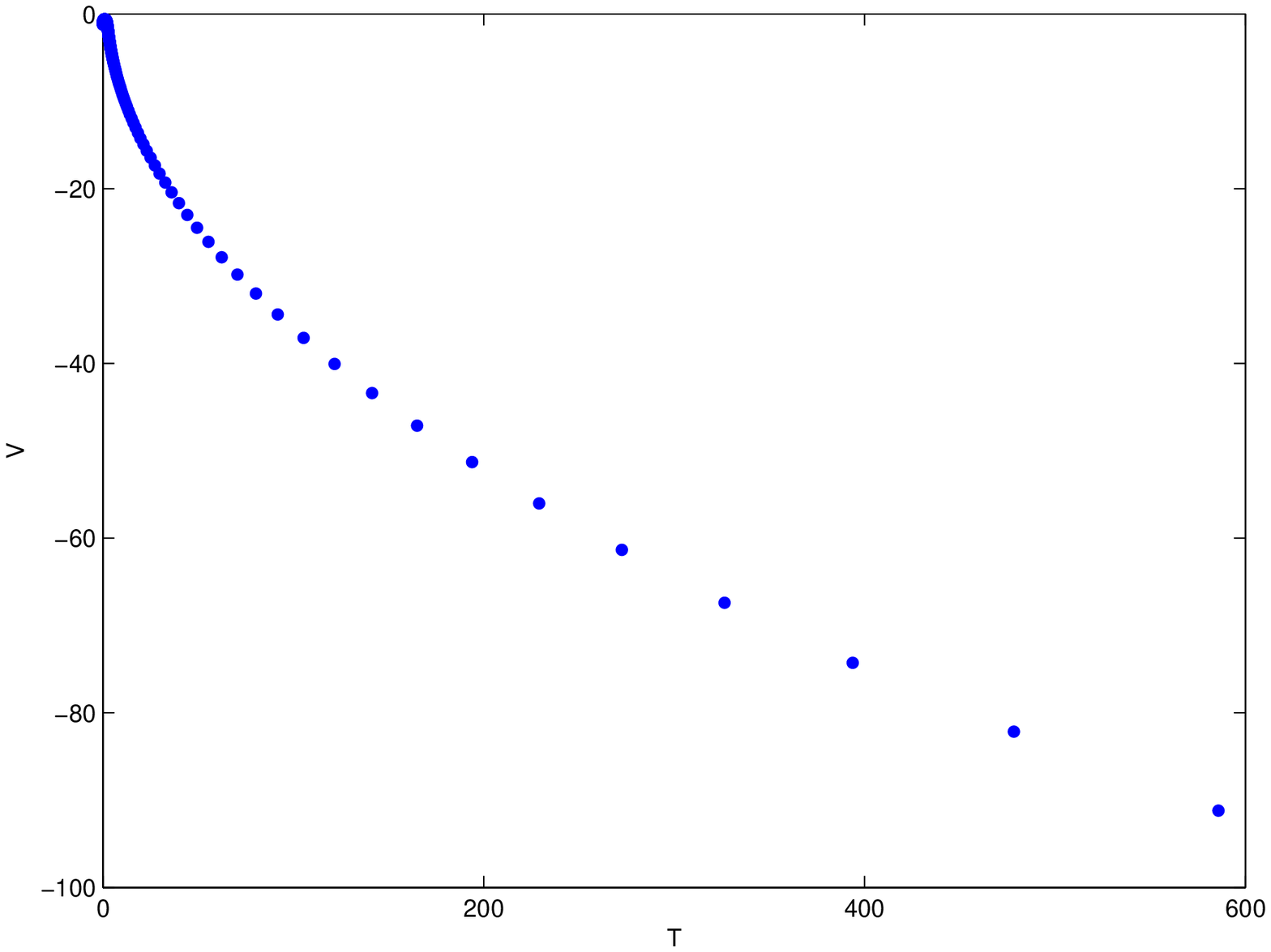}
\end{floatrow}
\caption{Behavior of the velocity $V$ versus $T$ given by the scaling transformation
(\ref{time-variable}) with $n = 1$ (left) and $n = 2$ (right).}
\label{fig-power-1}
\end{figure}

Figure \ref{fig-power-1} shows the dependence of $V$ versus the rescaled time variable $T$
for $n = 1$ (left) and $n = 2$ (right). It is obvious that the blow up occurs in
finite $T$ time if $n = 1$ and in infinite $T$ time if $n = 2$, which corroborates
well with the previous numerical data suggesting that $p < 0.5$.
This figure rules out the validity of the logarithmic law (\ref{loglaw}).
We have checked that the rescaled time variable $T$ for $n \geq 3$ also extends
to infinite times, similarly to the result for $n = 2$.

We note that the dependence of $V$ versus the original
time variable $t$ can be obtained by numerical integration of the integral in (\ref{time-variable}).
We have checked that both time evolutions of $V$ in $T$ with $n = 1$ and $n = 2$ recover the same behavior of
$V$ in $t$, which resembles the top left panel of Figure \ref{fig-num-1} except
times near the blow-up time, where the computational error becomes more significant.

Using the scaling transformation (\ref{time-variable}) with $n = 2$ in the case when
$T \to \infty$ as $t \to \infty$, we can define a more accurate procedure of detecting
the blow-up rate $p$ in the power law (\ref{pwrlaw}). First, we note
that if $|V(t)| = \mathcal{O}((t_0 - t)^{-p})$ as $t  \nearrow t_0$,
then $T = \mathcal{O}((t_0 - t)^{1 - 4p})$ as $t  \nearrow t_0$. Hence
$V(T) = \mathcal{O}(T^q)$ as $T \to \infty$ with $q := \frac{p}{4p-1}$.
Using now the linear regression in log-log variables for $T$ and $V$, we can
estimate the coefficient $q$, and then $p$. The following table shows
several computations of $q$ and $p$ for different initial and terminal
times $T$. All other parameters are fixed similarly to the previous numerical computations.

\begin{verbatim}
Starting time    Terminal time    Regression slope q    Blow-up rate p
     36.0943         723.3424           0.5345              0.4697
    121.7362         723.3424           0.5221              0.4797
    272.5828       78034.1670           0.5044              0.4956
   2393.6301       78034.1670           0.4997              0.5003
\end{verbatim}

The results of data fitting suggest that the power law \eqref{pwrlaw} gives a consistent estimation
of the blow-up rate $p$, with $p \approx 0.5$. Let us now discuss why the behavior
$|V(t)| = \mathcal{O}((t_0 - t)^{-1/2})$ as $t \nearrow t_0$ appears a generic feature of
smooth solutions of the boundary--value problem (\ref{pde2})--(\ref{bc2}).

Using equations (\ref{velocity1}) and (\ref{r5}), we obtain the dynamical equation
on $\beta(t) = u_x |_{x =0}$:
\begin{equation}
\label{contact}
\frac{d \beta}{d t} = -\frac{u_{xxxx} |_{x = 0}}{2 \beta(t)} - u_{xxxxx} |_{x = 0}, \quad t \geq 0.
\end{equation}
Let us now assume that there is $t_0 > 0$ such that
\begin{equation}
\label{assumed-rate}
\beta(t) \to 0, \quad u_{xxxx} |_{x = 0} \to a_4, \quad u_{xxxxx} |_{x = 0} \to a_5,
\quad \mbox{\rm as} \quad t \nearrow t_0,
\end{equation}
where $a_4 \neq 0$ and $|a_5| < \infty$. Solving the differential
equation (\ref{contact}) near the time $t = t_0$, we obtain
\begin{equation}
\label{asymptotic-rate}
\beta^2(t) = a_4(t_0 - t) + \mathcal{O}(t_0 - t)^{3/2}, \quad \Rightarrow \quad
V(t) = \sqrt{\frac{a_4}{t_0 - t}} + \mathcal{O}(1), \quad \mbox{\rm as} \quad t \nearrow t_0,
\end{equation}
under the constraint that $a_4 > 0$. The asymptotic rate (\ref{asymptotic-rate})
corresponds to the power law \eqref{pwrlaw} with $p = 0.5$.

\begin{figure}[H]
\begin{floatrow}
\includegraphics[height = 2.2 in, width = 3.2 in]{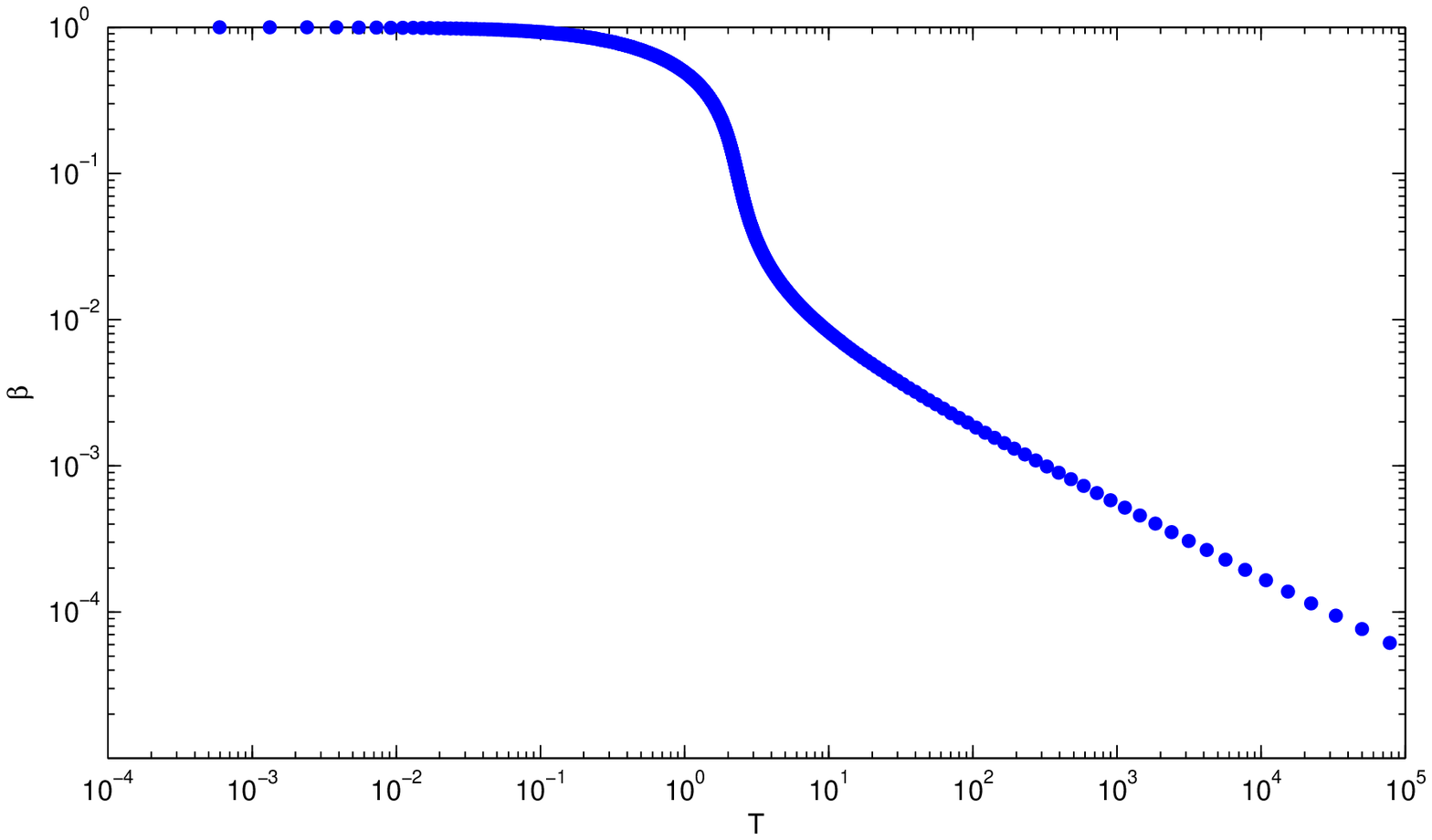}
\includegraphics[height = 2.2 in, width = 3.2 in]{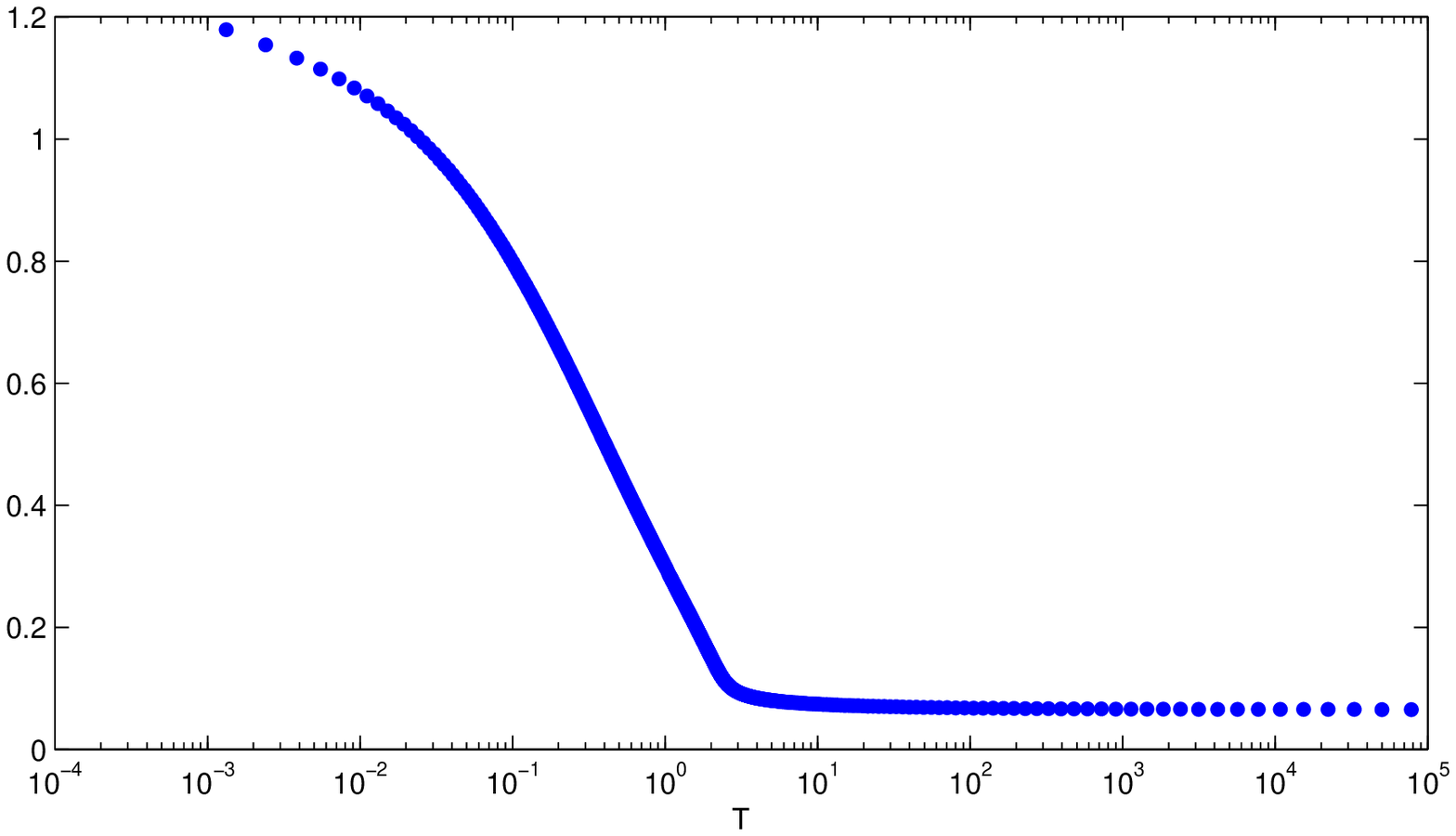}
\end{floatrow}
\caption{Behavior of $|u_x |_{x=0}|$ (left) and
$u_{xxxx} |_{x=0}$ (right) versus $T$ given by the scaling transformation
(\ref{time-variable}) with $n = 2$. The logarithmic scaling is used for $T$ and
$\beta$ variables.}
\label{fig-power-2}
\end{figure}

Figure \ref{fig-power-2} shows the behavior of absolute values of $u_x |_{x=0}$ (left) and
$u_{xxxx} |_{x=0}$ (right) versus the rescaled time variable $T$ given by
(\ref{time-variable}) with $n = 2$. 
We can see that the assumption $a_4 > 0$,
that is, $u_{xxxx} |_{x=0}$ is bounded away from zero near the blow-up time,
is justified numerically. We note that the time evolution in the rescaled time
variable $T$ allows us to identify this property better than the time evolution
in the original time variable $t$, which is shown on the bottom right panel of
Figure \ref{fig-num-1}. We have also checked from the linear regression
in log-log coordinates that $|\beta(t)| = \mathcal{O}((t_0 - t)^p)$ as $t \nearrow t_0$ with $p \approx 0.5$,
in consistency with the asymptotic rate (\ref{asymptotic-rate}).

\section{Conclusion}

We conclude from the numerical simulations of the boundary--value problem (\ref{pde2})--(\ref{bc2})
that, for any suitable initial condition in the two-parameter form \eqref{ic}, there always
exists a finite positive time $t_0$ such that $V(t) \rightarrow -\infty$ as $t  \nearrow t_0$,
although the blow-up time $t_0$ varies from different initial velocity $V(0)$.
With a large positive initial velocity $V(0)$, the solution tends to have a longer phase of slow motion
before it eventually blows up, whereas a negative initial velocity yields a much smaller value of the blow-up time.

The numerical results also suggest that the behavior of $V(t)$ near the blow-up time satisfies
the power law \eqref{pwrlaw}, with a blow-up rate $p \approx 0.5$. This numerical observation
corroborates a simple analytic theory for the blow-up of the velocity of contact lines
in the reduced model (\ref{pde})--(\ref{bc-pde}). Based on earlier numerical evidences
in \cite{Benilov}, a similar result should also hold for the nonlinear thin-film equation
(\ref{model}).

An open problem for further studies is to develop a more precise and computationally efficient numerical method
for solutions of the boundary--value problem (\ref{pde2})--(\ref{bc2}).
Because the model equation \eqref{pde2} is already a fourth-order differential equation,
we shall avoid using any numerical methods that involves higher-order central differences.
In addition, because of the unknown variable $V(t)$, it is difficult to use other higher-order
implicit methods to solve the system of differential equations after discretization. Thus,
the finite difference method has a limited accuracy. Therefore, a different approach
is needed, for instance, by using the collocation method involving the discrete Fourier transform.

\end{document}